\documentclass[11pt]{amsart}

\usepackage{amsrefs}

\newcommand {\R}  {\ensuremath{\mathbb{R}}}

\newcommand{\namedRef}[1]{%
\csname env:#1\endcsname\ \ref{#1}%
}

\def\envLRT#1#2{\expandafter\gdef\csname env:#1\endcsname{#2}}

\makeatletter
\newcommand{\namedLabel}[1]{\expandafter\xdef\csname env:#1\endcsname{\envName}%
\immediate\write\@auxout{\string\envLRT{#1}{\envName}}%
\label{#1}}
\newcommand{\NamedLabel}[1]{\expandafter\xdef\csname env:#1\endcsname{\envName}%
\immediate\write\@auxout{\string\envLRT{#1}{\envName}}%
\label{#1}}
\makeatother

\def\newType#1{\newtheorem{#1}[thm]{#1\gdef\envName{#1}}}
\newtheorem{thm}{Theorem}[section]
\newtheorem{Theorem}[thm]{Theorem\gdef\envName{Theorem}}
\newType{Lemma}
\newType{Corollary}
\newType{Proposition}
\newType{Addendum}
\newenvironment{Named Theorem}[1][xxx]{%
\gdef\envName{#1}\par\noindent{\bf #1.}\bgroup\sl }{\egroup}

\theoremstyle{definition}

\newtheorem{Example}[thm]{Example\gdef\envName{Example}}

\newType{Notation}
\newenvironment{Comment}[1][Remark]{\par\noindent{\bf #1.} }{}

\newbox\nbox
\gdef\Bar#1{\setbox\nbox = \hbox{$#1$}
\kern .08\wd\nbox{\overline{\hbox to .84\wd\nbox{\vphantom {$#1$}\hss}}}
\kern .08\wd\nbox\kern -\wd\nbox\box\nbox}

\def\DFT{$DFT$}

\begin{document}

\title{Impossible metric conditions on exotic $\R^4$'s}

\author{Laurence R. Taylor}
\address{Department of Mathematics, University of Notre Dame, Notre Dame, IN 46556}
\email{taylor.2@nd.edu}
\thanks{Partially supported by the N.S.F}

\begin{abstract}
There are many theorems in the differential geometry literature of the following sort.
{\par\noindent\narrower\sl
Let $M$ be a complete Riemannian manifold with some conditions
on various curvatures, diameters, volumes, etc.
Then $M$ is homotopy equivalent to a finite CW complex, or $M$ is
the interior of a compact, topological manifold with boundary.\par}

\noindent
At first glance it seems unlikely that such theorems have anything to say
about smooth manifolds homeomorphic to $\R^4$.
However, there is a common theme to all the proofs
which forbids the existence of such metrics 
on most (and possibly all) exotic $\R^4$'s.
\end{abstract}

\maketitle

\section{Definitions and the main result}
We say a smooth manifold $E$ homeomorphic to $\R^4$ satisfies the
\DFT\ condition (for De Michelis, Freedman and Taubes) provided that,
for every compact subset $K \subset E$, there exists an open neighborhood $U$ of $K$ 
such that
\begin{enumerate}
\item $K\subset U$
\item the closure of $U$, $\Bar U$, is homeomorphic to $D^4$
\item $U$ can be engulfed by itself rel $K$
\end{enumerate}
Precisely, condition (3) means that there exists a smooth, ambient, isotopy of $E$ 
from the identity to $\iota$ such that
$\iota\colon U\to E$ satisfies $\Bar U\subset\iota(U)$ and the isotopy is the identity on $K$.

As discussed in section 2, all \emph{exotic} $\R^4$'s (smoothings of $\R^4$ not diffeomorphic
to the standard smoothing) known to the author do not have the \DFT-property.

Differential geometry enters the picture via critical points of functions related to the distance function. 
Such functions are not necessarily smooth so the notion of critical point needs 
to be interpreted and there is a standard way to do this going back at least to Grove and Shiohama
\cite{Grove and Shiohama}.
The idea involves constructing vector fields and using differential geometry to get enough 
similarity to the gradient-like vector fields of Morse theory to prove results.
Here is the main observation of this note.
\begin{Theorem}\namedLabel{main result}
Let $E$ be a smooth, manifold homeomorphic to $\R^4$ with a proper 
Lipschitz function which has bounded critical values.
Then $E$ satisfies the \DFT-condition.
\end{Theorem}

\section{Remarks on the \DFT\ property}
Taubes \cite[Thm.~1.4, p.~366]{Taubes} proves that many exotics $\R^4$'s have the property 
that neighborhoods of large compact sets can not be embedded smoothly in exotic $\R^4$'s
with periodic ends. 
Hence these $\R^4$'s are not \DFT\ since $\iota$ can be used to construct a periodic end 
smoothing of $\R^4$ containing a neighborhood of any compact set $K$.

De Michelis and Freedman \cite{De Michelis and Freedman} state in the first line of the last paragraph
on page 220 that the family of $\R^4$'s they construct do not satisfy the \DFT-property
even if the $\iota$ is not required to be isotopic to the identity.

Gompf \cite{Gompf} argues that none of the smoothings of $\R^4$ in his menagerie have
the \DFT\ property. 
Certainly the universal $\R^4$ of Freedman and Taylor \cite{Freedman and Taylor} is not \DFT.

On the other hand, if the smooth Schoenflies conjecture is true but the smooth Poincar\'e 
conjecture is false, there will be exotic \DFT\ $\R^4$'s.

\section{Patterns of application}
Rather than produce a long list  of theorems, here is a meta-principle for generating theorems.
\begin{Named Theorem}[Meta-Principle]\namedLabel{Meta-Principle}
Take any theorem in differential geometry regarding the existence of complete metrics
with special properties. 
If the proof shows that a distance function or some other proper Lipschitz function 
has bounded critical values then such metrics do not exist on a non-\DFT$\,\R^4$.
\end{Named Theorem}

\begin{Comment}[Remark]%
The distance function from any point $p\in E$ is proper Lipschitz if the Riemannian metric on $E$
is complete.
\end{Comment}

Here are three examples.
\begin{Example}
Reading the paper of Lott and Shen \cite{Lott and Shen} gave the author the idea for this note.
The first paragraph on page 281 shows that an exotic $\R^4$ 
with lower quadratic curvature decay, quadratic volume growth and which does not collapse at
infinity is \DFT.
\end{Example}

\begin{Example}
An early finite-type theorem is by Abresch \cite{Abresch} which says that if the curvature
decays faster than quadratic on some complete Riemannian exotic $\R^4$, then it satisfies
the \DFT-property.
\end{Example}

\begin{Example}
An example involving mostly Ricci curvature and diameter is 
\cite[Thm.~B, p.~356]{Abresch and Gromoll}.
\end{Example}

\section{The proof of Theorem 1.1}
Let $E$ be any smoothing of $\R^4$.
Say that a flat topological embedding $e\colon S^3 \subset E$ is 
\emph{not a barrier to isotopy} provided there is a smooth vector field on $E$ 
with compact support so that if $\iota\colon E \to E$ is the isotopy at time $1$ generated
by the vector field, then
$\Bar{U} \subset \iota(U)$, where $U$ is the bounded component of $E - e(S^3)$.

Say that $E$ has \emph{no barrier to isotopy to $\infty$} if, for every compact set
$K\subset E$, there is a flat topological embedding $e\colon S^3 \subset E$ such that $K$ lies
in the bounded component of $E - e_K(S^3)$ and such that $e_K$ is not a barrier to isotopy .

\begin{Proposition}\namedLabel{no barrier to isotopy}
Let $E$ be a smooth, manifold homeomorphic to $\R^4$.
If $E$ has no barrier to isotopy to $\infty$ then  $E$ satisfies the \DFT-condition.
\end{Proposition}
\begin{proof}
Given $K$ compact, pick an $e_K\colon S^3 \to E$ with $K$ in the bounded component
of $E - e_K(S^3)$ such that $e_K$ is not a barrier to isotopy.
Let $\chi_1$ be the smooth vector field promised by the definition.
Let $U$ be the bounded component of $E - e_K(S^3)$ and notice $U$ satisfies (1) and (2).
Since $K \subset U$, there is another smooth vector field $\chi$ such that $\chi$
vanishes on $K$ and agrees with $\chi_1$ on a neighborhood of $e(S^3)$.
Let $I\colon E\times [0,1] \to E$ be the isotopy generated by the flow for $\chi$.
Since $\chi$ vanishes on $K$, the isotopy $I$ fixes $K$ and since $\chi$ agrees with
$\chi_1$ on $e(S^3)$, $I_1\bigl(e(S^3),1\bigr) \subset E - \Bar{U}$.
Hence $I$ is the isotopy required for (3).
\end{proof}

\begin{proof}[Proof of \namedRef{main result}]
Let $\rho\colon E \to [0,\infty)$ be a proper Lipschitz function.
Since the critical values are bounded by hypothesis, there is an $r_0$ such that 
for any critical point $x\in E$, $\rho(x) < r_0$.
Since $\rho$ is proper, $\rho^{-1}\bigl([0,r_0]\bigr)$ is compact.
Let $e\colon S^3 \to E$ be any flat embedding with $\rho^{-1}\bigl([0,r_0]\bigr)$ in
the bounded component of $E - e(S^3)$.
It will be shown that $e$ is not a barrier to isotopy, from which it follows that
$E$ has no barrier to isotopy to $\infty$, from which it follows that
$E$ satisfies the \DFT-property using \namedRef{no barrier to isotopy}.

Since $S^3$ is compact, $\rho\bigl(e(S^3)\bigr) \subset [r_1, r_2]$ with
$r_0<r_1$.
The idea is contained in the proofs  of \cite[Lemma 3.1, p. 108]{Greene} 
or \cite[Lemma 1.4, p.~2]{Cheeger}.
They start by constructing a vector field locally on $\rho^{-1}\bigl([r_1,r_2]\bigr)$ 
and patching it together using a smooth partition of unity. 
They observe that the conditions they need are open conditions so the field can be
taken to be smooth.
Then the proofs show that the resulting flow (or the flow for the negative of the constructed
field) moves $\rho^{-1}(r_1)$ out past $\rho^{-1}(r_2)$ and so $e(S^3)$ ends up in
$\rho^{-1}(r_2,\infty) = E - \rho^{-1}\bigl([0, r_2]\bigr)  \subset E - \Bar{U}$ since 
$\rho\bigl(\Bar{U}\bigr) \subset [0,r_2]$.
\end{proof}

\section{Concluding remarks}
Differential geometry can show that two smooth $4$-manifolds are diffeomorphic, 
avoiding the ``greater than or equal to $5$'' hypothesis of differential topology.
As examples, the Cartan-Hadamard Theorem \cite[Thm.~4.1, p.~221]{Sakai} 
and the Cheeger-Gromoll Soul Theorem \cite[Thm.~3.4, p.~215]{Sakai} both prove that
a $4$-manifold with very restrictive curvature conditions is 
\emph{diffeomorphic} to the standard $\R^4$.

The results presented here start with much weaker hypotheses than the
 Cartan-Hadamard or the Cheeger-Gromoll Soul Theorems, but the conclusions are also
 weaker.
 One question would be whether some of these theorems could be strengthened to 
 show that the manifold was diffeomorphic to the standard $\R^4$.
 A second question would be to use the \DFT\ property to produce interesting metrics, 
 perhaps metrics strong enough to prove that a \DFT\  $\R^4$ is standard.

\end{document}